\documentclass[letterpaper, 10pt, conference]{ieeeconf}
\usepackage{setspace, epsfig, psfrag, amssymb, amsmath, cite} \newfont{\smalll}{cmr8}

\def\IR{{\mbox{\hbox{I\hskip-.1em R}}}}

\def\IS{\hbox{I\hskip-.1em S}}

\def\IC{\hbox{C\hskip-
.5em\raise.5ex\hbox{$\scriptscriptstyle\mid$}}\ }
\def\Ic{\hbox{\smalll C\hskip-
.5em\raise.3ex\hbox{$\scriptscriptstyle\mid$}}\ }

\def\T={\buildrel {\scriptscriptstyle\triangle} \over =}

\def\sqr#1#2{{\vcenter{\vbox{\hrule height.#2pt\hbox{\vrule
width.#2pt height#1pt \kern#1pt\vrule width.#2pt}\hrule
height.#2pt}}}}

\def\block-diag{\mathop{\rm block{\scriptstyle -}diag}}

\def\pmbb#1{\setbox0=\hbox{#1}\raise 0.5ex\box0}

\def\norm#1{\|#1\|}

\newcommand{\bequ}{\begin{eqnarray}}
\newcommand{\eequ}{\end{eqnarray}}
\newcommand{\mT}{^\mathrm{T}}

\newcommand {\beq}      {\begin{equation}}
\newcommand {\eeq}      {\end{equation}}

\def\IR{{\mathbb R}}
\def\IC{{\mathbb C}}

\def\IS{{\mathbb S}}

\renewcommand\theenumi{\roman{enumi}}
\renewcommand\theenumii{\alph{enumii}}
\renewcommand\theenumiii{\arabic{enumiii}}
\renewcommand\theenumiv{\Alph{enumiv}}

\def\sc{s_{\rm c}}

\def\sc{s_{\rm c}}

\usepackage{fancyhdr}
\usepackage{epsfig,psfrag,amsmath}
\usepackage{srcltx,color}
\usepackage{theorem}
\usepackage{latexsym,amssymb,amsmath,amsfonts}
\usepackage{verbatim}
\usepackage{algorithm}
\usepackage{algorithmic}
\usepackage{cite}
\usepackage{subcaption}
\usepackage{caption}
\graphicspath{{figures/}}
\theoremstyle{plain} \theorembodyfont{\upshape}

\IEEEoverridecommandlockouts
\title{\Huge{ {Resilient Coordination of Networked Multiagent Systems Based on Distributed State Emulators}}
}
\author{Tansel Yucelen and Gerardo De La Torre
\thanks{
\newline\indent T. Yucelen is an Assistant Professor of the Mechanical and Aerospace Engineering Department and the Director of the Advanced Systems Research Laboratory at the Missouri University of Science and Technology, Rolla, MO 65409, USA (e-mail: {\tt\small yucelen@mst.edu}).  
\newline\indent G. De La Torre is a Graduate Research Assistant of the School of Aerospace Engineering at the Georgia Institute of Technology, Atlanta, GA 30332, USA (e-mail: {\tt\small glt3@gatech.edu}).
\newline\indent \hspace{-0.250cm} This research was supported by the University of Missouri Research Board.
}
}
  
\begin{document}  \maketitle \baselineskip 12.1pt
\maketitle

\definecolor{tRed}{RGB}{250,50,0}
\newcommand{\red}[1]{{\color{tRed}{{\textbf{#1}}}}}


\begin{abstract}

This note studies resilient coordination of networked multiagent systems in the presence of misbehaving agents, i.e., agents that are subject to adversaries modeled as exogenous disturbances. 
Apart from the existing relevant literature that make specific assumptions on the graph topology and/or the fraction of misbehaving agents, we present an adaptive control architecture based on distributed state emulators and show that the nominal networked multiagent system behavior can be retrieved even if all agents are misbehaving. 

\end{abstract}


\subsection{Introduction}

Although networked multiagent systems are envisioned to autonomously function in place of humans for repetitive, demanding, and often safety-critical missions the current level of controls technology is incapable of providing the needed usability and resiliency of multiagent systems. 
Because, the control algorithms of these systems are computed distributively without having a centralized entity monitoring the activity of agents, and hence, adversaries such as attacks to the communication network and/or failure of agent-wise components can easily result in system instability and prohibit the accomplishment of system-level objectives \cite{bullo2009distributed}.
The fragile nature of multiagent systems has triggered the development of detection and isolation algorithms during the last few years \cite{reference25,reference26,reference27,reference29,reference30,reference31}. 
For example, \cite{reference25} and \cite{reference26} make a specific assumption on the network connectivity (other than the standard assumption on the connectedness of networked agents) and \cite{reference27} requires that at most a fraction of any normal agent's neighbors to be adversaries, or misbehaving agents, for achieving resilient multiagent system behavior. 
Like \cite{reference25} and \cite{reference26}, a computationally expensive and not scalable algorithm is proposed in \cite{reference29} and \cite{reference30} based on input observers technique, where the effect of misbehaving agents on the overall multiagent system performance is also quantified, and an extension of this work is given in \cite{reference31} also focusing on the detection and isolation of misbehaving agents. 

This note develops an adaptive control architecture to ensure resilient coordination of networked multiagent systems in the presence of misbehaving agents. 
Specifically, we show that the nominal networked multiagent system behavior can be retrieved with the proposed methodology that utilizes a novel distributed state emulators. 
Apart from the existing relevant literature \cite{reference25,reference26,reference27,reference29,reference30,reference31} that make specific assumptions on the graph topology and/or the fraction of misbehaving agents, the proposed framework can achieve performance recovery on an arbitrary but connected communication topology and even if all agents are misbehaving. 


\subsection{Mathematical Preliminaries}

The notation used in this note is fairly standard.
Specifically, $\IR$ denotes the set of real numbers,
$\IR^{n \times m}$ denotes the set of $n \times m$ real matrices,
$\overline{\IS}^{n\times n}_+$ denotes the set of $n \times n$ symmetric nonnegative-definite real matrices,
$\mathrm{diag}(v)$ denotes a diagonal matrix with scalar (or matrix) entries given by $v$,
$\mathrm{I}_n$ denotes the $n\times n$ identity matrix,
$(\cdot)\mT$ denotes the transpose,
$\otimes$ denotes the Kronecker product, 
``$\triangleq$'' denotes the equality by definition, and 
$\mathbf{1}_{n}$ denotes a $n \times 1$ vector with 1 in all entries.
In addition, we write
$\norm{\cdot}$ for the euclidean norm,
$\norm{\cdot}_\textrm{F}$ for the Frobenius norm,
$\lambda_i(A)$ for the $i$-th eigenvalue of the matrix A,
$\textrm{spec}(A)$ for the spectrum of matrix A,
$\pi_0(A)$, $\pi_+(A)$, and $\pi_-(A)$ for the number of eigenvalues (counted with algebraic multiplicities) of A having zero, positive, and negative real parts, respectively, 
$\mathrm{det}(A)$ for the determinant of A,
$\mathrm{max}(\cdot)$ for the maximum, and
$\mathrm{min}(\cdot)$ for the minimum.

We now introduce several results that are necessary for the development of the main result of this note. 

\textbf{Lemma 1} \cite{blockm}\textbf{.} Consider the matrix given by
\bequ
M=\left[\begin{array}{cc}A & B \\ 0 & D\end{array}\right]. 
\eequ 
Then, the determinant of $M$ satisfies
\begin{equation} 
\textrm{det}(M)=\textrm{det}(A)\textrm{det}(D).
\end{equation}

\textbf{Lemma 2} [\citen{yuc02}]\textbf{.} Suppose $Z(\lambda)=A\lambda^2+B\lambda+C$ denotes the quadratic matrix polynomial where $A\in\mathbb{R}^{n \times n}$ and $C\in\mathbb{R}^{n \times n}$, and $A$ is nonsingular.
If $B\in\mathbb{R}^{n\times n}$ is positive-definite, then $\pi_+(Z)=\pi_-(A)+\pi_-(C),$ $\pi_-(Z)=\pi_+(A)+\pi_+(C),$ and $\pi_0(Z)=\pi_0(C),$ where $\pi_+(Z)+\pi_-(Z)+\pi_0(Z)=2n.$

We next recall some of the basic notions from graph theory and networked multiagent systems  [\citen{pre1},\citen{pre2}]. 
Specifically, graphs are broadly used in networked multiagent systems to encode interactions between a group of agents. 
An \emph{undirected} graph $\mathcal{G}$ is defined by a set $\mathcal{V}_\mathcal{G}=\{1,\dots,n\}$ of \emph{nodes} and a set $\mathcal{E}_\mathcal{G}\subset\mathcal{V}_\mathcal{G}\times\mathcal{V}_\mathcal{G}$ of \emph{edges}.
If the \emph{unordered} pair $(i,j)\in\mathcal{E}_\mathcal{G}$, then nodes $i$ are $j$ are \emph{neighbors} and the neighboring relation is indicated with $i\sim j$. 
The set of neighbors of node $i$ is denoted by $\mathcal{N}_\mathcal{G}(i)=\{j|(i\sim j)\in\mathcal{E}(\mathcal{G})\}$.
The \emph{degree} of a node is given by the number of its neighbors. 
Letting $d_i$ be the degree of node $i$, then the \emph{degree} matrix of a graph $\mathcal{G}$, ${\Delta(\mathcal{G})}\in\mathbb{R}^{n\times n}$, is given by ${\Delta(\mathcal{G})}\triangleq \mathrm{diag}(d), d=[d_1,\dots,d_n]^\textrm{T}$.  
The \emph{adjacency} matrix of a graph $\mathcal{G}$, $\mathcal{A(G)}\in\mathbb{R}^{n \times n}$, is given by
\begin{equation}
[ {A(\mathcal{G})} ]_{ij}\triangleq\left\{
     \begin{array}{cl}
       1, & \mathrm{if}\ (i,j)\in\mathcal{E}_\mathcal{G},\\
       0, & \mathrm{otherwise}.
     \end{array}
   \right.
\end{equation}
The \emph{Laplacian} matrix of a graph, $\mathcal{L(G)}\in\bar{\mathbb{S}}_+^{n \times n}$, plays a central role in many graph theoretic treatments of networked multiagent systems is given by $\mathcal{L(G)}\triangleq{\Delta(\mathcal{G})}-\mathcal{A(G)}$, 
where the spectrum of the Laplacian for a connected, undirected graph can be ordered as $0=\lambda_1(\mathcal{L(G)})<\lambda_2(\mathcal{L(G)})\leq\cdots\leq\lambda_n(\mathcal{L(G)})$, with $\mathbf{1}_n$ as the eigenvector corresponding to the zero eigenvalue $\lambda_1(\mathcal{L(G)})$ and $\mathcal{L(G)}\mathbf{1}_n=\mathbf{0}_n$ holds.

Networked multiagent systems can be modeled by a graph $\mathcal{G}$, where nodes and edges, represent agents and interagent information exchange links, respectively.
In particular, let $x_i(t)\in\mathbb{R}^{N}$ denote the state of agent $i$ at time $t\geq0$, whose dynamics are described by
\begin{equation} \label{single_agent}
\dot{x}_i(t)=u_i(t),\quad t\geq0, \quad x_i(t_0)=x_{i0}, \quad i\in\mathcal{V}(\mathcal{G}),
\end{equation}
with $u_i(t)\in\mathbb{R}^N$ being the control input of agent $i$. 
We consider agents having dynamics of the form given by (\ref{single_agent}) to focus the main result of this note.
In addition, we focus on the consensus problem without loss of much generality when presenting the main contribution of this note. 
In particular, if agent $i$ is allowed to access the relative state information with respect to its neighbors, a solution to the standard consensus problem can be given by
\begin{equation}\label{input}
u_i(t)=-\sum_{i\sim j}(x_i(t)-x_j(t)), \quad i\in\mathcal{V}(\mathcal{G}), 
\end{equation}
for a connected and undirected graph (throughout this note, we assume that the graph $\mathcal{G}$ is connected and undirected). 
The networked multiagent system given by (\ref{single_agent}) and (\ref{input}) can now be described in the form given by
\vspace{0cm}
\begin{equation}\label{complex}
\dot{x}(t)=-(\mathcal{L(G)}\otimes\mathrm{I}_N)x(t), \quad t\geq0, \quad x(t_0)=x_0,
\end{equation}
where $x(t)=[x^\textrm{T}_1(t),\cdots,x^\textrm{T}_n(t)]^\textrm{T}$ denotes the aggregated state vector.
For ease of exposition, we consider the case of $N=1$.  
However, all results presented in this note can be trivially extended to the general case. 
Finally, considering ($\ref{complex}$), we note that 
\bequ 
x(t)\rightarrow[(\mathbf{1}_n\mathbf{1}^\textrm{T}_n/n)\otimes\mathrm{I}_N]x_0 \ \ \textrm{as} \ \ t\rightarrow\infty, 
\eequ 
since the graph $\mathcal{G}$ is assumed to be connected and undirected.
That is, the networked multiagent system is said to reach a consensus since $x_1=x_2=\cdots=x_n$.


\subsection{Resilient Coordination Based on State Emulators and Adaptive Control}

This section introduces the proposed adaptive control approach based on state emulators to enable resilient coordination of networked multiagent systems in the presence of misbehaving agents.
The agent dynamics given by (\ref{single_agent}) are augmented to incorporate the effect of these misbehaviors as
\begin{equation} \label{single_agent_with_constant_disturbance}
\dot{x}_i(t)=u_i(t)+w_i,\ t\geq0,\ x_i(t_0)=x_{i0}, \ i\in\mathcal{V}(\mathcal{G}), 
\end{equation}
where $w_i\in\mathbb{R}$ is an unknown disturbance applied to agent $i$. 
Notice that we represent adversaries as disturbances similar to \cite{yuc01}.  
Specifically, we say that an agent is misbehaving if there exists a time such that $w_i\neq0$, $ i\in\mathcal{V}(\mathcal{G})$. 
It should be mentioned here that we only consider the case of constant exogenous disturbances for the ease of exposition (using the results of Section VI of \cite{yuc01}, our following results can be easily extended to the case of time-varying disturbances).

In order to mitigate the effect of these exogenous disturbances, the nominal consensus protocol given by (\ref{input}) is modified as
\begin{equation}\label{input_constant_disturbance}
u_i(t)=-\sum_{i\sim j}(x_i(t)-x_j(t))+\hat{w}_i, \quad i\in\mathcal{V}(\mathcal{G}), 
\end{equation}
where $\hat{w}_i$ is an adaptive control signal, which estimates the disturbance of agent $i$, and is updated as
\begin{equation}\label{hatw_update}
\dot{\hat{w}}_i(t)=\alpha(x_i(t)-\hat{x}_i(t)), \quad t\geq0, \quad \hat{w}_i(t_0)=\hat{w}_{i0},
\end{equation}  
where $\alpha>0$, $\hat{w}_{i0}=0$, and $\hat{x}_i(t)$ is a state emulator given by  
\begin{equation}\label{observer}
\dot{\hat{x}}_i(t)=-\textrm{d}_i\hat{x}_i(t)+\sum_{i\sim j}x_j(t) ,\quad t\geq0, \quad \hat{x}_i(0)=\hat{x}_{i0}.
\end{equation}
The undisturbed response of the system is captured by $\hat{x}_i(t)$.
However, notice that indirect disturbances from neighboring agents still affect the emulator system.   
Using (\ref{single_agent_with_constant_disturbance}) and (\ref{observer}), the dynamics of the emulator state estimate error are given by 
\begin{equation}\label{observer_error}
\dot{\tilde{x}}_i(t)=-d_i\tilde{x}_i(t)-\tilde{w}_i(t), \quad t\geq0, \quad \tilde{x}_i(t_0)=\tilde{x}_{i0},
\end{equation}
where $\tilde{x}(t)\triangleq x_i(t)-\hat{x}_i(t)$ and $\tilde{w}(t)\triangleq w_i-\hat{w}_i(t)$.
In addition, the exogenous disturbance estimate error dynamics are given by 
\begin{equation}\label{estimate_error}
\dot{\tilde{w}}_i(t)=-\alpha \tilde{x}_i(t),\quad t\geq0, \quad \tilde{w}_i(t_0)=\tilde{w}_{i0}.
\end{equation}
Now, the networked multiagent system can be described in a compact form as
\begin{eqnarray}\label{agg_sys}
\dot{\hat{x}}(t)&=&-\Delta(\mathcal{G})\hat{x}(t)+\mathcal{A(G)}{x}(t),\quad \hat{x}(t_0)=\hat{x}_{0},\label{agg_sys1} \ \ \ \\
\dot{\tilde{x}}(t)&=&-{\Delta(\mathcal{G})}\tilde{x}(t)-\tilde{w}(t), \quad \tilde{x}(t_0)=\tilde{x}_{0},\label{agg_sys2} \\
\dot{\tilde{w}}(t)&=&\alpha\tilde{x}(t), \quad \tilde{w}(t_0)=\tilde{W}_{0},\label{agg_sys3}
\end{eqnarray}
where $\hat{x}(t)=[\hat{x}_1(t), \dots, \hat{x}_n(t)]^\textrm{T}\in\mathbb{R}^n$, 
$\tilde{x}(t)=[\tilde{x}_1(t), \dots, \tilde{x}_n(t)]^\textrm{T}\in\mathbb{R}^n$,   
$\tilde{w}(t)=[\tilde{w}_1(t), \dots, \tilde{w}_n(t)]^\textrm{T}\in\mathbb{R}^n$,
denote the aggregated emulator state, emulator estimate error, and disturbance estimate error, respectively.   
Furthermore, (\ref{agg_sys1}) can be equivalently written as 
\begin{equation} \label{new_refer}
\dot{\hat{x}}(t)=-\mathcal{L}(\mathcal{G})\hat{x}(t)+\mathcal{A}(\mathcal{G})\tilde{x}(t), \ \ t\geq0, \ \  \hat{x}(t_0)=\hat{x}_{0}.
\end{equation}
Next, we consider the state transformation given by 
\begin{equation} \label{state_trans}
\hat{y}(t)=T\hat{x}(t)=[\hat{z}_{12}(t), \hat{z}_{13}(t), \dots, \hat{c}_\mathcal{G}(t)]^\textrm{T}
\end{equation}
where $\hat{z}_{1i}(t)=\hat{x}_1(t)-\hat{x}_i(t)$ and $\hat{c}_\mathcal{G}(t)=\sum_{i\in\mathcal{V}(\mathcal{G})}\hat{x}_i(t)$. 
Under this state transformation and using (\ref{new_refer}), it follows that 
\begin{equation} \label{trans_new_reference}
\dot{\hat{y}}(t)=-T\mathcal{L}T^{-1}\hat{y}(t)+T\mathcal{A}\tilde{x}(t), \ t\geq0, \ \hat{y}(t_0)=\hat{y}_{0}. 
\end{equation}
Furthermore, since the graph is undirected and connected, (\ref{trans_new_reference}) can be partitioned as
\begin{eqnarray} 
\hspace{-0.4cm}\dot{\hat{z}}_1(t)&=&A_1\hat{z}_1(t)+A_2\tilde{x}(t),\ t\geq0, \ \hat{z}_1(t_0)=\hat{z}_{10},     \label{nonsemistable_dynamics}\\
\hspace{-0.4cm}\dot{\hat{c}}_\mathcal{G}(t)&=&\sum_{i\in\mathcal{V}(\mathcal{G})}\textrm{d}_i\tilde{x}_i(t),\ \hat{c}_\mathcal{G}(t_0)=\hat{c}_{\mathcal{G}0},\label{centroid_dynamics}
\end{eqnarray}
where $z_1=[z_{12},z_{13},\dots,z_{1n}]^\textrm{T},$ $A_1\in\mathbb{R}^{n-1\times n-1}$ and $A_2\in\mathbb{R}^{n-1\times n-1}$ are the matrices obtained by removing the $n^{\textrm{th}}$ row and column from $T\mathcal{L}(\mathcal{G})T^{-1}$ and the matrix obtain by removing the $n^{\textrm{th}}$ row from $T\mathcal{A}(\mathcal{G})$, respectively.

The estimate update given by (\ref{hatw_update}) requires each agent to have access to its state, $x_i(t)$. 
This requirement is not needed for the standard consensus protocol (\ref{input}).
However, as shown in the following result, the modified consensus protocol (\ref{input_constant_disturbance}) allows the network multiagent system to achieve consensus even in the presence of disturbances.   
Specifically, the following result shows that the modified consensus protocol given by (\ref{input_constant_disturbance}) results in consensus despite the presence of disturbances.
Furthermore, the estimates of the disturbances converge to the actual exogenous disturbances as $t\rightarrow\infty$. 

\textbf{Theorem.} Consider the network multiagent system given by (\ref{nonsemistable_dynamics}), (\ref{agg_sys2}), and (\ref{agg_sys3}).
Then, the solution $(\hat{z}_1(t),\tilde{x}(t),\tilde{w}(t))$ is exponentially stable for all $(\hat{z}_{10}, \tilde{x}_{0}, \tilde{w}_0)\in\mathbb{R}^{n-1}\times\mathbb{R}^{n}\times\mathbb{R}^{n}$. 

\textbf{Proof.} Note that the system is equivalently described by $\dot{\xi}(t)=M\xi(t),\quad t\geq0, \quad {\xi}(t_0)={\xi}_{0}$, 
where 
\begin{eqnarray}\label{aug_agg_sys_matrix}
M=\left[\begin{array}{c|cc}
  A_1 & A_2 & 0_{n-1\times n} \\
  \hline
  0_{n\times n-1} & -{\Delta(\mathcal{G})} & -\textrm{I}_n \\
0_{n\times n-1} & \alpha\textrm{I}_n & 0_{n\times n} \\
 \end{array}\right],
\end{eqnarray}
and $\xi(t)=[\hat{z}_1(t)^\textrm{T}, \tilde{x}(t)^\textrm{T}, \tilde{w}(t)^\textrm{T}]^\textrm{T}\in\mathbb{R}^{3n-1}$.
Using Lemma 1 the spectrum of $M$ is described as  
\begin{eqnarray}
\textrm{spec}(M)=\textrm{spec}(A_1)\cup\textrm{spec}(\left[\begin{array}{cc}-{\Delta(\mathcal{G})} & -\textrm{I}_n \\ \alpha\textrm{I}_n & 0 \\ \end{array}\right]).
\end{eqnarray}
Since the graph is assumed to be connected the spectrum of the Laplacian is described as  
\begin{eqnarray}
\textrm{spec}(-\mathcal{L(G)})&=&\{0\}\cup\{\lambda_2(-\mathcal{L(G)}),\dots, \lambda_n(-\mathcal{L(G)})\}\nonumber\\
&=&\{0\}\cup\textrm{spec}(A_1)
\end{eqnarray}
where $\lambda_i(-\mathcal{L(G)})<0, \forall i\in\{2,\dots,n\}$.
Furthermore, note that the characteristic polynomial of $\left[\begin{array}{cc}-{\Delta(\mathcal{G})} & -\textrm{I}_n \\ \alpha\textrm{I}_n & 0 \\ \end{array}\right]$ is given as $Z(\lambda)=\lambda^2\textrm{I}_n+\lambda \Delta(\mathcal{G}) +\alpha \textrm{I}_n$. 
Therefore, it can be concluded from Lemma 2 that $\pi_+(Z)=0$, $\pi_0(Z)=0$ and $\pi_-(Z)=2n$. 
Thus, $\lambda_i({M})<0, \forall i\in\{1,\dots,3n-1\}$. 
Therefore, the system is exponentially stable for all initial conditions and $t\geq0$.
\hfill $\blacksquare$ 

Notice from (\ref{agg_sys1}) and (\ref{centroid_dynamics}) that $\tilde{x}(t)$ acts as a vanishing perturbation to an ideal consensus equation.
Furthermore, if $\|\tilde{x}(t)\|_2$ is sufficiently small, then agents not only achieve consensus but the agreement point stays close to the original centroid of the system -- the point that would have been reached in the undisturbed case.
That is, the effect of the disturbances on the overall system performance can be related to $\|\tilde{x}(t)\|_2.$
In addition, $\alpha$ has a direct effect on the bound of $\|\tilde{x}(t)\|_2$.
To see this, consider the energy function $E(\tilde{x}(t),\tilde{w}(t))=\frac{1}{2}\tilde{x}^\textrm{T}\tilde{x}+\frac{1}{2\alpha}\tilde{w}^\textrm{T}\tilde{w}$. 
Taking the time derivative yields $\dot{E}(\tilde{x}(t),\tilde{w}(t))=-\tilde{x}^\textrm{T}\Delta(\mathcal{G})\tilde{x}\leq0$. 
Therefore, $E(\tilde{x}(t),\tilde{w}(t))\leq E(\tilde{x}(0),\tilde{w}(0))=\frac{\|\tilde{w}\|^2_2}{2\alpha}$ and $\|\tilde{x}\|_2\leq\frac{\|\tilde{w}\|_2}{\sqrt{\alpha}}$. 
As $\alpha$ is increased, the magnitude of the vanishing perturbation term $\|\tilde{x}\|_2$ becomes smaller.
Meaning, during transient-time the state emulator system (\ref{agg_sys1}) and the emulator centroid, $\hat{c}(t)$, are effected less by disturbances. 
Finally, note that $\hat{c}_\mathcal{G}(t)=\int_0^t\sum_{i\in\mathcal{V}(\mathcal{G})}\textrm{d}_i\tilde{x}_i(t) \ \textrm{d}t + \hat{c}_{0}$ remains bounded since $\tilde{x}(t)$ exponentially converges to $0$. 


\subsection{Concluding Remarks}
Control algorithms of networked multiagent systems are generally computed distributively without having a centralized entity monitoring the activity of agents; and therefore, adversaries such as attacks to the communication network and/or failure of agent-wise components can easily result in system instability and prohibit the accomplishment of system-level objectives. 
Motivation from this standpoint, we proposed a new adaptive control approach based on distributed state emulators to guarantee a desired system-level performance in the presence of misbehaving agents. 


\bibliographystyle{IEEEtran}
\bibliography{Reference}

\begin{thebibliography}{10}
\providecommand{\url}[1]{#1}
\csname url@samestyle\endcsname
\providecommand{\newblock}{\relax}
\providecommand{\bibinfo}[2]{#2}
\providecommand{\BIBentrySTDinterwordspacing}{\spaceskip=0pt\relax}
\providecommand{\BIBentryALTinterwordstretchfactor}{4}
\providecommand{\BIBentryALTinterwordspacing}{\spaceskip=\fontdimen2\font plus
\BIBentryALTinterwordstretchfactor\fontdimen3\font minus
  \fontdimen4\font\relax}
\providecommand{\BIBforeignlanguage}[2]{{%
\expandafter\ifx\csname l@#1\endcsname\relax
\typeout{** WARNING: IEEEtran.bst: No hyphenation pattern has been}%
\typeout{** loaded for the language `#1'. Using the pattern for}%
\typeout{** the default language instead.}%
\else
\language=\csname l@#1\endcsname
\fi
#2}}
\providecommand{\BIBdecl}{\relax}
\BIBdecl

\bibitem{bullo2009distributed}
F.~Bullo, J.~Cortes, and S.~Martinez, \emph{Distributed control of robotic
  networks: a mathematical approach to motion coordination algorithms}.\hskip
  1em plus 0.5em minus 0.4em\relax Princeton University Press, 2009.

\bibitem{reference25}
S.~Sundaram and C.~N. Hadjicostis, ``Distributed function calculation via
  linear iterations in the presence of malicious agents: Attacking the
  network,'' in \emph{American Control Conference, 2008}.\hskip 1em plus 0.5em
  minus 0.4em\relax IEEE, 2008, pp. 1350--1355.

\bibitem{reference26}
------, ``Distributed function calculation via linear iterations in the
  presence of malicious agents: Overcoming malicious behavior,'' in
  \emph{American Control Conference, 2008}.\hskip 1em plus 0.5em minus
  0.4em\relax IEEE, 2008, pp. 1356--1361.

\bibitem{reference27}
H.~J. LeBlanc, H.~Zhang, S.~Sundaram, and X.~Koutsoukos, ``Resilient
  continuous-time consensus in fractional robust networks,'' in \emph{American
  Control Conference (ACC), 2013}.\hskip 1em plus 0.5em minus 0.4em\relax IEEE,
  2013, pp. 1237--1242.

\bibitem{reference29}
F.~Pasqualetti, A.~Bicchi, and F.~Bullo, ``Distributed intrusion detection for
  secure consensus computations,'' in \emph{Decision and Control, 2007 46th
  IEEE Conference on}.\hskip 1em plus 0.5em minus 0.4em\relax IEEE, 2007, pp.
  5594--5599.

\bibitem{reference30}
------, ``Consensus computation in unreliable networks: A system theoretic
  approach,'' \emph{Automatic Control, IEEE Transactions on}, vol.~57, no.~1,
  pp. 90--104, 2012.

\bibitem{reference31}
I.~Shames, A.~M. Teixeira, H.~Sandberg, and K.~H. Johansson, ``Distributed
  fault detection for interconnected second-order systems,'' \emph{Automatica},
  vol.~47, no.~12, pp. 2757--2764, 2011.

\bibitem{blockm}
J.~R. Silvester, ``Determinants of block matrices,'' \emph{The Mathematical
  Gazette}, vol.~84, no. 501, pp. 460--467, 2000.

\bibitem{yuc02}
B.~Bilir and C.~Chicone, ``A generalization of the inertia theorem for
  quadratic matrix polynomials,'' \emph{Linear Algebra and its Applications},
  vol. 280, no. 2–3, pp. 229 -- 240, 1998.

\bibitem{pre1}
M.~Mesbahi and M.~Egerstedt, \emph{Graph theoretic methods in multiagent
  networks}.\hskip 1em plus 0.5em minus 0.4em\relax Princeton University Press,
  2010.

\bibitem{pre2}
A.~Rahmani, M.~Ji, M.~Mesbahi, and M.~Egerstedt, ``Controllability of
  multi-agent systems from a graph-theoretic perspective,'' \emph{SIAM Journal
  on Control and Optimization}, vol.~48, no.~1, pp. 162--186, 2009.

\bibitem{yuc01}
T.~Yucelen and M.~Egerstedt, ``Control of multiagent systems under persistent
  disturbances,'' in \emph{American Control Conference}, 2012, pp. 5264--5269.

\end{thebibliography}


\begin{thebibliography}{1}
\providecommand{\url}[1]{#1}
\csname url@samestyle\endcsname
\providecommand{\newblock}{\relax}
\providecommand{\bibinfo}[2]{#2}
\providecommand{\BIBentrySTDinterwordspacing}{\spaceskip=0pt\relax}
\providecommand{\BIBentryALTinterwordstretchfactor}{4}
\providecommand{\BIBentryALTinterwordspacing}{\spaceskip=\fontdimen2\font plus
\BIBentryALTinterwordstretchfactor\fontdimen3\font minus
  \fontdimen4\font\relax}
\providecommand{\BIBforeignlanguage}[2]{{%
\expandafter\ifx\csname l@#1\endcsname\relax
\typeout{** WARNING: IEEEtran.bst: No hyphenation pattern has been}%
\typeout{** loaded for the language `#1'. Using the pattern for}%
\typeout{** the default language instead.}%
\else
\language=\csname l@#1\endcsname
\fi
#2}}
\providecommand{\BIBdecl}{\relax}
\BIBdecl

\bibitem{pre1}
M.~Mesbahi and M.~Egerstedt, \emph{Graph theoretic methods in multiagent
  networks}.\hskip 1em plus 0.5em minus 0.4em\relax Princeton University Press,
  2010.

\bibitem{pre2}
A.~Rahmani, M.~Ji, M.~Mesbahi, and M.~Egerstedt, ``Controllability of
  multi-agent systems from a graph-theoretic perspective,'' \emph{SIAM Journal
  on Control and Optimization}, vol.~48, no.~1, pp. 162--186, 2009.

\bibitem{blockm}
J.~R. Silvester, ``Determinants of block matrices,'' \emph{The Mathematical
  Gazette}, vol.~84, no. 501, pp. 460--467, 2000.

\bibitem{yuc01}
T.~Yucelen and M.~Egerstedt, ``Control of multiagent systems under persistent
  disturbances,'' in \emph{American Control Conference}, 2012, pp. 5264--5269.

\bibitem{yuc02}
B.~Bilir and C.~Chicone, ``A generalization of the inertia theorem for
  quadratic matrix polynomials,'' \emph{Linear Algebra and its Applications},
  vol. 280, no. 2–3, pp. 229 -- 240, 1998.

\bibitem{Khalil}
H.~K. Khalil, \emph{Nonlinear Systems}.\hskip 1em plus 0.5em minus 0.4em\relax
  Prentice Hall, 1996.

\bibitem{Freq-limit}
T.~Yucelen, G.~D.~L. Torre, and E.~Johnson, ``Improving transient performance
  of adaptive control architectures using frequency-limited system error
  dynamics,'' \emph{arXiv:1309.6693}, vol.~48, no.~1, pp. 162--186, 2009.

\bibitem{r:bernstein}
D.~S. Bernstein, \emph{Matrix Mathematics: Theory, Facts, and Formulas},
  2nd~ed.\hskip 1em plus 0.5em minus 0.4em\relax Princeton, NJ: Princeton
  University Press, 2009.

\bibitem{Khalil:2002}
H.~K. Khalil, \emph{Nonlinear Systems}.\hskip 1em plus 0.5em minus 0.4em\relax
  Upper Saddle River, NJ: Prentice-Hall, 2002.

\end{thebibliography}
\end{document}